\documentclass[11pt,a4paper]{article} 

\usepackage{amsmath,amsthm,amssymb,latexsym} 
\usepackage{amsfonts} 
\usepackage[american]{babel} 
\usepackage{mathrsfs} 
\usepackage[dvips]{graphicx} 
 
\newtheorem{teo}{Theorem}[section] 
\newtheorem{prop}[teo]{Proposition} 
\newtheorem{lemma}[teo]{Lemma}%[section] 
\newcommand{\R}{\mathbb{R}}

\newcommand{\bd}{\textrm{bd}} 
\newcommand{\cnp}{c_{n,p}} 
\newcommand{\conv}{\textrm{conv}} 
 
\newcommand{\dpkm}{\kmax+_p(-\kmax)} 
\newcommand{\ee}{\mathbf{e}} 
\newcommand{\enu}{e_{n+1}} 
\newcommand{\hk}{h_{K}} 
\newcommand{\ii}{\mathscr{I}} 
\newcommand{\intt}{\textrm{int}} 
\newcommand{\kdue}{\mathscr{K}^2} 
\newcommand{\kmax}{\overline{K}} 
\newcommand{\kn}{\mathscr{K}^n} 
\newcommand{\kt}{\tilde{K}} 
\newcommand{\la}{\lambda} 
\newcommand{\lt}{\tilde{L}} 
\newcommand{\pp}{\mathscr{P}} 
\newcommand{\ttt}{\mathscr{T}} 
\newcommand{\vn}{V_n} 
 
\newcommand{\mink}{Minkowski} 
 
\newcommand{\shs}{shadow system} 
\newcommand{\lps}{linear parameter system} 
\newcommand{\cm}{continuos movement} 
\newcommand{\pcm}{parallel chord movement} 
\newcommand{\res}{Rogers \& Shephard} 
\begin{document} 
%% 
%%%%%%%%%%%%%%%%%%%%%%%%%%%%%%%%%%%%%%%%%%%%%%%%%%%%%%%%%%%%%%%%%%%%%%%%%%%%%%%%%%%%%%%%%%%%%%% 
%%%%%%%%%%%%%%%%%%%%%%%%%%%%%%%%%%%%%%%%%%%%%%%%%%%%%%%%%%%%%%%%%%%%%%%%%%%%%%%%%%%%%%%%%%%%%%% 
%% 
\begin{title}{ 
A sharp {\res} type inequality for the $p$-difference body of planar 
convex bodies. }\end{title} 
 
\author{C. Bianchini and A. Colesanti} 
\date{} 
\maketitle 
 
\begin{abstract} 
We prove a sharp {\res} type inequality for the $p$-difference body 
of a convex body in the two-dimensional case, for every $p\ge 1$. 
\end{abstract} 
 
\bigskip 
 
%% 
%%------------------------------------------------------------------------------------------- 
\section{Introduction} 
%%------------------------------------------------------------------------------------------- 
%% 
A \emph{convex body} is a non-empty convex compact subset of $\R^n$; 
let us indicate the set of convex bodies in $\R^n$ with $\kn$. To 
each convex body $K$ we can associate in a biunique way its 
\emph{support function} $h_K$: 
$$ 
\hk(u)=\sup\{ \langle x,u\rangle \;|\;x\in K \},\qquad\textrm{ for all } u\in\R^n, 
$$ 
where $\langle\cdot,\cdot\rangle$ is the standard scalar product. 
The support function is a fundamental tool since the main properties of the body can be deduced from it. 
 
One of the most interesting aspects of convex geometry, i.e. the 
theory of convex bodies, are geometric inequalities. An important 
family of inequalities are those leading to estimate the volume of a 
special body associated with a convex body (for example the 
difference body or the  reflection body) in terms of the volume of 
the body itself. 
 
A remarkable inequality of this type is the classical {\res} inequality (see \cite{RS1}) which asserts that for all $K\in\kn$ 
\begin{equation}\label{resin} 
\vn\Big(K+(-K)\Big)\le \binom{2n}{n} \vn(K), 
\end{equation} 
and equality holds if and only if $K$ is a simplex. 
Here $\vn(K)$ denotes the $n$-dimensional volume of $K$ (i.e. the $n$-dimensional Lebesgue measure). 
The body $K+(-K)$ is called difference body of $K$ and it is the {\mink} sum of $K$ and its reflected body with respect to the origin, $-K$. We recall more generally that the {\mink} sum of $K$ and $L\in\kn$ is 
$$ 
K+L=\{ z\in \R^n\;|\; z=k+l,\quad k\in K,\,l\in L \}. 
$$ 
 
Another inequality due to {\res} (\cite{RS2}) concerns the convex hull (here denoted by $\conv$) of $K$ and $-K$, under the assumption that the origin $o$ belongs to $K$: 
\begin{equation}\label{infresin} 
\vn\Big(\conv\big(K\cup(-K)\big)\Big)\le 2^n \vn(K), 
\end{equation} 
where equality holds if and only if $K$ is a simplex with one vertex at the origin. 
 
In \cite{F} Firey introduced a new operation for convex bodies, 
called $p$-sum, which depends on the parameter $p\ge 1$ and extends 
the {\mink} sum. An account on the theory of convex bodies based on 
the $p$-sum, the so called \emph{Brunn-{\mink}-Firey theory}, can be 
found in the papers \cite{L1}, \cite{L2} by Lutwak. Let us fix 
$K,L\in\kn$ both containing the origin; the $p$-sum of $K$ and $L$, 
$K+_pL$, is defined by its support function in the following way: 
$$ 
h_{K+_pL}(u)=\Big( \hk^p(u)+h_L^p(u) \Big)^\frac 1p,\qquad u\in\R^n. 
$$ 
This definition admits a natural extension to the case $p=\infty$: 
$$ 
h_{K+_{\infty}L}(u)=\lim_{p\to\infty}h_{K+_pL}(u)=\max\{ h_K(u),h_L(u) \},\; u\in \R^n. 
$$ 
 
Note that the extremal values $p=1$ and $p=\infty$, correspond to 
the {\mink} sum and the convex hull of the union respectively. 
Indeed one has 
$$ 
h_{K+_1L}(u)=h_K(u)+h_L(u)=h_{K+L}(u), 
$$ 
and 
$$ 
h_{K+_{\infty}L}(u)=\max\{ h_K(u),h_L(u) \}=h_{\conv(K\cup L)}. 
$$ 
 
As proved by Firey \cite{F}, the $p$-sum is monotone with respect to the parameter $p$: for all $K,L\in\kn$ such that $o\in K,L$, if $p\le q$ then 
$$ 
K+_qL\subseteq K+_pL. 
$$ 
This implies that for all $p\ge 1$, 
$$ 
\conv(K\cup L)\subseteq K+_pL \subseteq K+L. 
$$ 
Another simple inclusion is 
$$ 
K+_pL\subseteq 2^{\frac 1p}\conv(K\cup L). 
$$ 
 
In particular choosing $L=-K$ and using inequalities (\ref{resin}) and (\ref{infresin}), we have: 
$$ 
\vn\Big( K+_p(-K) \Big)\le \min\Big\{ \binom{2n}{n},\; 
2^{n\frac{(1+p)}{p}} \Big\}\;\vn(K). 
$$ 
A natural problem is then to find the best constant $c=c_{n,p}$, depending on $n$ and $p$, such that 
\begin{equation}\label{presn} 
\vn\Big( K+_p(-K) \Big)\le \;c_{n,p}\;\vn(K),\qquad\textrm{for all }K\in\kn,\,o\in K. 
\end{equation} 
 
In this paper we solve this problem in the planar case $n=2$, for every $p\ge 1$. 
\begin{teo}\label{presteo} 
For every $p\ge 1$ there exists a constant $c_p$ such that 
\begin{equation}\label{pres} 
V_2\Big(K+_p(-K)\Big)\le c_p\,V_2(K), 
\end{equation} 
for all $K\in\kdue$. 
In particular if $K$ is a triangle with one vertex at the origin, then equality holds. 
\end{teo} 
An explicit expression of $c_p$ will be presented in Section \ref{secproof}. 
 
We will show the $p$-{\res} inequality (\ref{pres}) as a consequence 
of a theorem about the $p$-sum of the so called {\pcm}s of convex 
bodies. 
 
A parallel chord movement is a special one-parameter family of 
convex bodies which can be seen as continuous deformations of a 
fixed convex body. More precisely, fix $K\in\kn$ and a direction 
$v\in\R^n$ which is the direction of the movement. We move each 
chord of $K$ parallel to $v$ in that direction with a certain speed 
and we consider the union of these chords as the time parameter 
varies. If the speed function is suitably chosen, namely if the 
union of the chords is convex for all values of the parameter, then 
the family of the resulting convex bodies is a {\pcm}. 
 
Parallel chord movements are special cases of a wider class of 
movements of convex bodies introduced by Rogers and Shephard in 
\cite{RS3}, which have been recently applied in the proof of several 
inequalities in convex geometry (see, for examples, \cite{CCG}, 
\cite{CG1}-\cite{CG5}, \cite{MR}). 
 
The importance of these movements is due principally to the 
behaviour  of several geometric functionals with respect to the 
parameter of the movement. Indeed many of them, and the volume is 
the main example, are convex function of the time parameter of the 
movement. 
 
In particular in this paper we prove that if $K_t$ is a {\pcm}, then 
the volume of its $p$-difference body $\vn\Big(K_t+_p(-K_t)\Big)$ is 
a convex function of $t$, for all $p\ge 1$. This result, together 
with a thecnique used in \cite{CCG}, leads to the proof of Theorem 
\ref{presteo}. As noted in \cite{CCG} this thecnique is successful 
only in the planar case, so our method can not be used to prove 
inequality (\ref{presn}) in the general case $n\ge 2$. 
 
The paper is organized as follows. In Section \ref{secsslps} we introduce several kinds of movements of convex bodies and we show some of their properties. Next to basic results we present a theorem about the $p$-sum of a particular type of movements. In section \ref{secproof} we prove Theorem \ref{presteo} as an application of the results concerning movements of convex bodies. 
%% 
%%------------------------------------------------------------------------------------------- 
\section{Shadow systems and \lps s}\label{secsslps} 
%%------------------------------------------------------------------------------------------- 
%% 
A \emph{{\shs}} is a family of $n$-dimensional convex bodies $\{K(u)\}$ obtained as the projection of a fixed convex body $\kt\subseteq \R^{n+1}$ onto the hyperplane $\{\enu^{\perp}\}$, which we identify with $\R^n$, along the direction $\enu +u$. Here $u$ varies in $\{\enu^{\perp}\}$. The {\shs} is said to be originated from the $(n+1)$-dimensional body $\kt$. 
 
A \emph{{\lps}} is a family of convex bodies $\{K_t\}$ that can be written in the form 
\begin{equation}\label{lpsdef} 
K_t= \conv\{ x_i+\la_itv\;:\;i\in I \}, \qquad {t\in\ii}; 
\end{equation} 
where $I$ is an arbitrary index set, $\{x_i\}_{i\in I}$ and $\{\la_i\}_{i\in I}$ are bounded subsets of $\R^n$ and of $\R$ respectively, $\ii$ is an interval of $\R$ and $v\in\R^n$ is the direction of the {\lps}. 
 
Linear parameter systems are \shs s in which $u$ lies on a line, indeed we have the following result. 
\begin{prop}\label{lpseshs} 
$\{K_t\}_{t\in\ii}$ is a {\lps} in $\R^n$ if and only if there exists a convex body $\kt$ in $\R^{n+1}$ such that for every $t\in\ii$, $K_t$ is the projection of $\kt$ onto the hyperplane $\{\enu^{\perp}\}$ along the direction $\enu-tv$. 
\end{prop} 
The idea to view \lps s as projections of higher dimensional convex 
bodies is contained in the original papers by Rogers and Shephard 
(\cite{RS3}, \cite{S}) and was largely used by Campi and Gronchi 
(\cite{CG1}-\cite{CG5}). For the sake of completeness here we 
present the proof of Proposition \ref{lpseshs}. 
\begin{proof} 
Let $K_t$ be of the form (\ref{lpsdef}) and  let us define the body $\kt$ as follows: 
$$ 
\kt=conv\Big( \{ x_i+\la_i\enu\ :\ i\in I\ \} \Big). 
$$ 
For all $t\in\ii$ let us call $L_t$ the projection of $\kt$ onto $\{\enu^{\perp}\}$ along $\enu-tv$. For all $y\in L_t$ there exists $z\in \kt$ such that $ y=z-\langle z,\enu\rangle\big(\enu-tv\big)$. Furthemore there exist $a_i\in \enu^{\perp}$, $\lambda_i\in\R$, and $\sigma_i\ge 0$, $i=1,...,n+1$, such that $\sum_{i=1}^{n+1}\sigma_i=1$ and 
$$ 
z=\sum_{i=1}^{n+1}\sigma_i(a_i+\lambda_i\enu). 
$$ 
Therefore 
$$ 
y=\sum_{i=1}^{n+1}\sigma_i(a_i+\lambda_itv). 
$$ 
This implies that $L_t$ is contained in $K_t$. To prove the reverse inclusion one can observe that the previous implications are true in both directions. 
 
Conversely, let $\kt$ be any $(n+1)$-dimensional convex body and fix $t\in\ii$; its projection onto $\{\enu^{\perp}\}$ along $\enu-tv$ is the set: 
$$ 
L_t=\{\enu^{\perp}\}\cap\{z\in\R^{n+1}\ |\ z=x+t\xi,\ x\in\kt,\ t\in\R \}. 
$$ 
This is equivalent to write 
$$ 
L_t= \big\{ x-\langle x,\enu\rangle\enu+t\langle -x,\enu\rangle w\ :\ x\in\kt \big\}, 
$$ 
and, by the convexity of $\kt$, $\{L_t\}_{t\in\ii}$  is a {\lps} as defined in (\ref{lpsdef}). 
\end{proof} 
From the previous proof it follows that the body $\kt$ which generates a {\lps} of the form (\ref{lpsdef}) can be explicitely written as: 
\begin{equation}\label{ktilde} 
\kt=\conv\{x_i+\la_i \enu\;:\;i\in I\}. 
\end{equation} 
Campi and Gronchi showed in \cite{CG4} the following formula which relates the support functions of $K_t$ and $\kt$: 
\begin{equation}\label{funzsup} 
h_{K_t}(u)=h_{\kt}(u+t\langle u,v\rangle\enu),\qquad u\in\R^n,\; t\in\ii. 
\end{equation} 
 
We can give a cinematic interpretation of a {\lps} viewing the 
numbers  $\la_i$ as the speeds of the points $x_i$ along the 
direction $v$ and $t$ as the time parameter. 
 
If the index set $I$ is a convex body $K\in\kn$ and the speed is a 
function of the point, then the  {\lps} is called \emph{{\cm}}: 
$$ 
K_t=\conv \{ x+\alpha(x)tv\;:\; x\in K \},\qquad t\in\ii, 
$$ 
where $\alpha(\cdot)$ is a bounded function on $K$. 
 
Assume that the speed function is constant on each chord parallel to 
$v$, i.e. $\alpha(x)=\beta(x|v^{\perp})$ where $x|v^{\perp}$ is the 
projection of $x$ onto $\{v^{\perp}\}$ and $\beta$ is a function 
defined on the ortoghonal projection of $K$ onto $\{v^{\perp}\}$. 
Moreover, if $\beta$ is such that convexity is preserved for any 
$t$, namely 
$$ 
\{ x+\beta(x|v^{\perp})tv\;:\; x\in K \}=\conv \{ x+\beta(x|v^{\perp})tv\;:\; x\in K \}, 
$$ 
then the {\cm} is called \emph{\pcm}. 
 
% A {\pcm} such that the speed function is constant on each chord parallel to the direction $v$ is called \emph{{\pcm}}. 
% A {\pcm} along the direction $v$ is 
% $$K_t= \{ x+\beta(x|v^{\perp})tv\;:\; x\in K \},\qquad t\in\ii;$$ 
% where $x|v^{\perp}$ is the projection of $x$ onto $\{v^{\perp}\}$. 
 
In other words a {\pcm} is obtained assigning to each chord parallel 
to the direction $v$ a speed vector $\beta(x|v^{\perp})v$ and 
considering for each fixed time $t$ the union of these chords. Such 
union has to be convex. Notice that if $\{K_t\}_{t\in\ii}$ is a 
{\pcm}, then the volume of $K_t$ is independent of $t$. 
 
The following theorem is due to Rogers and Shephard (see \cite{RS3}) 
and it is one of the main motivations for the use of \lps s in the 
theory of convex bodies. %% 
\begin{teo}\label{convvol} 
The volume $\vn(K_t)$ of a {\lps} is a convex function of the parameter $t$. 
\end{teo} 

In \cite{CG4} it is proved that the {\mink} sum of \lps s is a {\lps}. Here we extend this result to the $p$-sum. This fact is one of the main ingredients in the proof of the $p$-{\res} inequality. 
\begin{teo}\label{lpspsum} 
Let $\{K_t\}_{t\in\ii}$ and $\{L_t\}_{t\in\ii}$ be \lps s along the direction $v$ and let $p\ge 1$, then $\{K_t+_pL_t\}_{t\in\ii}$ is also a {\lps} along the direction $v$. 
\end{teo} 
The proof is a straightforward consequence of Proposition \ref{lpseshs} and  the following lemma. 
\begin{lemma} 
Let $\{K_t\}_{t\in\ii}$ and $\{L_t\}_{t\in\ii}$ be \lps s along the same direction $v$ and let $\kt$ and $\lt$ be the $(n+1)$-dimensional convex bodies which generate $K_t$ and $L_t$ respectively, defined as in (\ref{ktilde}). Hence for all $t\in\ii$, $K_t+_pL_t$ is the projection of $\kt+_p\lt$ onto the hyperplane $\{\enu^{\perp}\}$ along the direction $\enu-tv$. 
\end{lemma} 
\begin{proof} 
Using (\ref{funzsup}) one has: 
\begin{eqnarray*} 
h_{\kt+_p\tilde{L}}^p\big(u+t\langle u,v\rangle \enu\big)  &=&  h_{\kt}^p\big(u+t\langle u,v 
\rangle \enu\big) + h_{\tilde{L}}^p\big(u+t\langle u,v 
\rangle \enu\big)\\ 
               &=&  h_{K_t}^p(u)+h_{L_t}^p(u)= h_{K_t+_pL_t}^p(u). 
\end{eqnarray*} 
This implies that $K_t+_pL_t$ is the projection of the body $\kt+_p\lt$ onto the hyperplane $\{\enu^{\perp}\}$ along the direction $\enu-tv$, which means, by Proposition (\ref{lpseshs}), that $K_t+_pL_t$ is a {\lps} along $v$. 
\end{proof} 
%% 
%%------------------------------------------------------------------------------------------- 
\section{The proof of the $p$-{\res} inequality}\label{secproof} 
%%------------------------------------------------------------------------------------------- 
%% 
Let us call $\kn_0$ the set of convex bodies with non-empty interior and containing the origin and let us consider the functional $F_p$ defined on $\kn_0$: 
$$ 
F_p(K)=\frac{ \vn\Big( K+_p(-K) \Big) }{ \vn(K) }. 
$$ 
 
It is clear that the best constant $\cnp$ such that (\ref{presn}) holds is the supremum of $F_p$ in $\kn_0$. 
 
% The functional is continuos on $\kn_0$ and non-singular linear transformations invariant so  using John's theorem \cite{J} and the Blaschke selection's theorem, the existence of a  maximizer for $F_p$ follows. 
 
We will use \lps s to find a maximum for the functional $F_p$ in the planar-case. 
The starting point is the next proposition which follows from Theorem \ref{convvol} and Theorem \ref{lpspsum}. 
\begin{prop}\label{convfp} 
If $K_t$ is any {\pcm} such that $K_t\in\kn_0$ for all $t\in\ii$, then $F_p(K_t)$ is a convex function of the parameter $t$. 
\end{prop} 

In \cite{CCG} the following fact is proved: if $P$ is a planar convex polygon with $m$ vertices, $m>3$, then there exists a {\pcm}  $\{P_t\}_{t\in[t_0,t_1]}$, with $t_0<0<t_1$, such that $P=P_0$ and $P_{t_0}$ and $P_{t_1}$ have at most $(m-1)$ vertices. 
By Proposition \ref{convfp} it follows that: 
$$ 
F_p(P)\le \max\{ F_p(P_{t_0}), F_p(P_{t_1}) \}. 
$$ 
Using recursively this fact we deduce that: 
$$ 
\sup_{\pp}F_p=\sup_{\ttt}F_p, 
$$ 
where $\pp=\{K\in\kdue_0\;|\; K \textrm{ is a polygon }\}$ and $\ttt=\{K\in\kdue_0\;|\; K \textrm{ is  triangle }\}$. Moreover, by the continuity of $F_p(\cdot)$ and a standard density argument, one has: 
$$ 
\sup_{\kdue_0}F_p=\sup_{\ttt}F_p. 
$$ 
 
In particular we are going to show that triangles with one vertex at the origin are maximizers for $F_p$. 
In order to do this, let $T\in\ttt$ and assume $o\in\intt(T)$ ($\intt$ denotes the interior). Then there exists a {\pcm} (whose elements are translates of $T$), $\{T_t\}_{t\in[t_0,t_1]}$ with $t_0<0<t_1$, such that $T_0=T$ and $o\in\bd(T_{t_0})$, $o\in\bd(T_{t_1})$, ($\bd$ denotes the boundary). 
Similarly, if $o\in\bd(T)$, then there exists a {\pcm} containing $T$, whose endpoints are triangles with one vertex at $o$. 
Using again Proposition \ref{convfp}, we have proved that 
$$ 
\sup_{\kdue_0}F_p=\sup_{\ttt_0}F_p, 
$$ 
where $\ttt_0$  the set of triangles with one vertex at the origin. 
 
Note that $F_p$ is invariant under non-singular linear transformations. 
This implies that $F_p$ is constant on $\ttt_0$. 
 
This argument proves the following result. 
\begin{teo} 
If $T$ is a triangle in $\kdue_0$ with one vertex at the origin, then $T$ is a maximizer for $F_p$. 
\end{teo} 
 
To compute the best constant $c_{2,p}$, we can choose as a maximizer the triangle with vertices at the origin, at $(1,0)$ and $(0,1)$; let us indicate it with $\kmax$. 
Namely 
$$ 
c_{2,p}=\frac{ V_2\Big( \kmax+_p(-\kmax) \Big) }{ V_2(\kmax) }. 
$$ 
 
Then to express the value of $c_{2,p}$ it is necessary to know how the $p$-difference body $\kmax+_p(-\kmax)$ looks like. Here we use  the parametrization of the boundary of a convex body in terms of its support function (see \cite{sch} Corollary 1.7.3). 
 
The support function of $\kmax+_p(-\kmax)$ is: 
$$ 
h_{\kmax+_p(-\kmax)}(w)= 
\begin{cases} 
\cos\theta & \qquad\textrm{ if }0\le\theta<\frac{\pi}{4},\\ 
\sin\theta & \qquad\textrm{ if }\frac{\pi}{4}\le\theta<\frac{\pi}{2},\\ 
\big( \sin^p\theta+(-\cos\theta)^p \big)^{\frac1p} & \qquad\textrm{ if }\frac{\pi}{2}\le\theta\le\pi, 
\end{cases} 
$$ 
where $w=\ee^{i\theta}\in S^1$. Furthemore, by the symmetry of $\kmax+_p(-\kmax)$, 
$$ 
h_{\dpkm}(\ee^{i\theta})=h_{\dpkm}(\ee^{i(\theta-\pi)}), 
$$ 
for all $\pi\le\theta\le 2\pi$. 
Then a parametrization for the boundary of $\kmax+_p(-\kmax)$, for $1<p<+\infty$, is $\zeta(\theta)=\Big( x(\theta), y(\theta) \Big)$, where 
$$ 
x(\theta)= 
\begin{cases} 
 1-\frac{2}{\pi}\theta & \qquad\textrm{ for }\theta\in[0,\frac{\pi}{2}],\\ 
 -\big(\sin^p\theta+(-\cos\theta)^p\big)^{\frac{1-p}{p}}(-\cos\theta)^{p-1} & \qquad\textrm{ for }\theta\in(\frac{\pi}{2},\pi);\\ 
\end{cases} 
$$ 
$$ 
y(\theta)= 
\begin{cases} 
 \frac{2}{\pi}\theta & \qquad\textrm{ for }\theta\in[0,\frac{\pi}{2}],\\ 
 \big(\sin^p\theta+(-\cos\theta)^p\big)^{\frac{1-p}{p}}\sin^{p-1}\theta & \qquad\textrm{ for }\theta\in(\frac{\pi}{2},\pi);\\ 
\end{cases} 
$$ 
and the remaining part of the boundary can be found using the symmetry of the body. 
 
A picture can perhaps better show the geometry of the body. In the following one $\kmax+_p(-\kmax)$ is represented for the the values $1$, $1.5$, $2$, $15$, $\infty$ of the parameter $p$. 
\begin{center} 
\includegraphics[height=4cm]{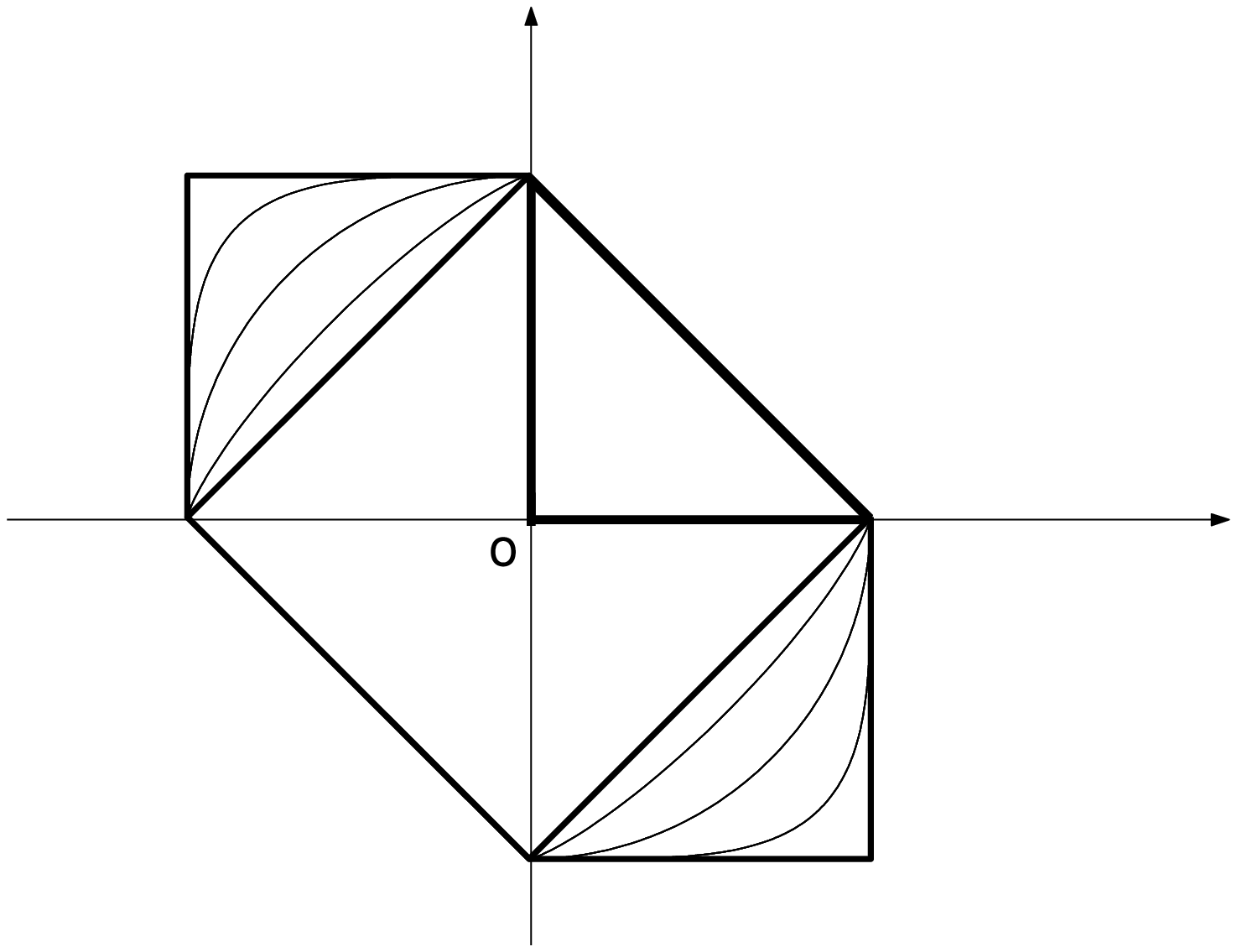} 
\end{center} 
 
Using the above parametrization and Gauss-Green's formulas we can express the area of $\dpkm$ and then the value of the best constant $c_{2,p}$: 
\begin{equation*}\label{c2p2} 
c_{2,p}=2\Bigg( 1+(p-1)\int_0^{\frac{\pi}{2}} \frac {\sin^{p-2}t\cos^{p-2}t} {\Big(\sin^pt+\cos^pt\Big)^{2\frac{(p-1)}{p}}} \, dt\Bigg),\quad1<p<+\infty. 
\end{equation*} 
% 
%% 
%% 
%%%%%%%%%%%%%%%%%%%%%%%%%%%%%%%%%%%%%%%%%%%%%%%%%%%%%%%%%%%%%%%%%%%%%%%%%%%%%%%%%%%% 
%% 
%% 
 
%% 
%% 
%%%%%%%%%%%%%%%%%%%%%%%%%%%%%%%%%%%%%%%%%%%%%%%%%%%%%%%%%%%%%%%%%%%%%%%%%%%%%%%%%%%%%%%%%%%%%%% 
\noindent {\em Authors' addresses} 
 
\noindent Dipartimento di Matematica ``U. Dini'', viale Morgagni 
67/A, 50134 Firenze, Italia. 
\medskip 
 
Chiara Bianchini: chiara.bianchini@math.unifi.it 
 
Andrea Colesanti: andrea.colesanti@math.unifi.it 
%%%%%%%%%%%%%%%%%%%%%%%%%%%%%%%%%%%%%%%%%%%%%%%%%%%%%%%%%%%%%%%%%%%%%%%%%%%%%%%%%%%%%%%%%%%%%%% 
%% 

\begin{thebibliography}{99} 
%% 
\addcontentsline{toc}{chapter}{Bibliografia} 
%% 
% 
\bibitem{CCG} \textsc{S.~Campi, A.~Colesanti, P.~Gronchi}, \textsl{A note on Sylvester's problem for random polytopes in a convex body}, Rend. Ist. Mat. Univ. Trieste \textbf{31} (1999), 79-94. 
% 
\bibitem{CG1} \textsc{S.~Campi, P.~Gronchi}, \textsl{The $L^p$-Busemann-Petty centroid inequality}, Adv. Math. \textbf{167} (2002), 128-141. 
% 
\bibitem{CG2} \textsc{S.~Campi, P.~Gronchi}, \textsl{On the reverse $L^p$-Busemann-Petty centroid inequality}, Mathematika \textbf{49} (2002), 1-11. 
% 
%\bibitem{CG3} \textsc{S.~Campi, P.~Gronchi}, \textsl{Extremal convex sets for Sylvester-Busemann type functionals}, Appl. Anal. \textbf{85} (2006), 129-141. 
% 
\bibitem{CG4} \textsc{S.~Campi, P.~Gronchi}, \textsl{On volume product inequalities for convex sets}, Proc. Amer. Math. Soc. \textbf{134}, 8 (2006), 2393-2402. 
% 
\bibitem{CG5} \textsc{S.~Campi, P.~Gronchi}, \textsl{Volume inequalities for $L_p$-zonotopes} (to appear on Mathematika). 
% 
\bibitem{F} \textsc{W.M.J.~Firey}, \textsl{p-Means of convex bodies}, Math. Scand. \textbf{10} 
(1962), 17-24. 
% 
\bibitem{J} \textsc{F.~John}, \textsl{Extremum problems with inequalities as subsidiary conditions}, Courant Anniversary Volume, Interscience, New York,  1948, 187-204. 
% 
\bibitem{L1} \textsc{E.~Lutwak}, \textsl{The Brunn-{\mink}-Firey theory I. Mixed volumes and the {\mink} problem}, J. Differential Geom. \textbf{38} (1993), 131-150. 
% 
\bibitem{L2} \textsc{E.~Lutwak}, \textsl{The Brunn-{\mink}-Firey theory II. Affine and geominimal surface area}, Adv. Math. \textbf{118}, 2 (1996), 244-294. 
% 
\bibitem{MR} \textsc{M.~Meyer, S.~Reisner}, \textsl{Shadow systems and volumes
    of polar convex bodies}, preprint (2006), arXiv:math.MG/0606305. 
% 
\bibitem{RS1} \textsc{C.A.~Rogers, G.C.~Shephard}, \textsl{The difference body of a convex body}, Arch. Math. \textbf{8} (1957), 220-233. 
% 
\bibitem{RS2} \textsc{C.A.~Rogers, G.C.~Shephard}, \textsl{Convex bodies associated with a given convex body}, J. Lond. Math. Soc. \textbf{33} (1958), 270-281. 
% 
\bibitem{RS3} \textsc{C.A.~Rogers, G.C.~Shephard}, \textsl{Some extremal problems for convex bodies}, Mathematika \textbf{5} (1958), 93-102. 
% 
\bibitem{sch} \textsc{R.~Schneider}, \textsl{Convex bodies: the Brunn-{\mink} theory}, Encyclopedia of Mathematics and its applications \textbf{44}, Cambridge University Press, Cambridge, 19993. 
% 
\bibitem{S} \textsc{G.C.~Shephard}, \textsl{Shadow system of convex sets}, Israel J. Math. \textbf{2} (1964), 229-236. 
% 
\end{thebibliography}
\end{document}